\documentclass[12pt]{amsproc}
\usepackage{amsfonts}
\usepackage{graphicx}

\def\N{\mathbb {N}}

\def\Z{\mathbb {Z}}
\def\R{\mathbb {R}}
\def\C{\mathbb {C}}
\def\Cbar{\overline{\C}}
\def\Cstar{\C^*}

\def\disk{\mathbb {D}}

\def\diskbar{\ovl{\disk}}

\def\ovl{\overline}
\def\phi{\varphi}
\def\eps{\varepsilon}
\def\theta{\vartheta}
\def\Re{\mbox{\rm Re}}
\def\Im{\mbox{\rm Im}}

\def\sm{\setminus}

\newtheorem{theorem}{Theorem} 

\newtheorem{lemma}[theorem]{Lemma}

\def\proof{\par\medskip\noindent {\sc Proof. }}
\def\proofof #1 {\par\medskip\noindent {\sc Proof of #1. }}

\def\sketchof #1 {\par\medskip\noindent {\sc Sketch of proof of #1. }}
\def\Box{\framebox[10pt]{\rule{0pt}{3pt}}}
\def\nix{\rule{0pt}{2pt}}
\def\qed{\qedd\par\medskip\noindent}
\def\qedd{\nix\nolinebreak\hfill\hfill\nolinebreak$\Box$}

\def\lineclear
    {\rule{0pt}{0pt}\nopagebreak\par\nopagebreak\noindent}

\def\reminder #1 {{\sf #1}}
\def\hide #1 {}

\def\ray{g_{\s}}
\def\rayp{g_{\s'}}
\def\raypp{g_{\s''}}
\def\rayki{g_{\s^{k,i}}}
\def\u{{\tt u}}
\def\uu{{\underline\u}}
\def\s{{\underline s}}

\def\su{\s}
\def\Sym{\mathcal{S}}

\begin{document}

\title[Iterated Cosine Maps]{The Dynamical Fine Structure of Iterated
Cosine Maps and a Dimension Paradox}
\author{Dierk Schleicher}
\address{School of Engineering and Science, IUB: International
University Bremen, Postfach 750 561, D-28725 Bremen, Germany}
\email{dierk@iu-bremen.de}
\date{\today}

\def\thefootnote{}
\footnotetext{{\em Date:} \today}
\def\thefootnote{\arabic{footnote}}

\begin{abstract}
We discuss in detail the dynamics of maps $z\mapsto ae^z+be^{-z}$ for
which both critical orbits are strictly preperiodic. The points which
converge to $\infty$ under iteration contain a set $R$ consisting of
uncountably many curves called ``rays'', each connecting $\infty$ to a
well-defined ``landing point'' in $\C$, so that every point in $\C$ is
either on a unique ray or the landing point of finitely many rays.

The key features of this paper are the following two: (1) this is the
first example of a transcendental dynamical system where the Julia set
is all of $\C$ and the dynamics is described in detail using symbolic
dynamics; and (2) we get the strongest possible version (in the plane) of
the ``dimension paradox'': the set $R$ of rays has Hausdorff dimension
$1$, and each point in $\C\sm R$ is connected to $\infty$ by one or more
disjoint rays in $R$; as a complement of a $1$-dimensional set, $\C\sm R$
has of course Hausdorff dimension $2$ and full Lebesgue measure. 
\end{abstract}

\maketitle

%\tableofcontents

\section{Introduction}
\label{SecIntro}

The dynamics of iterated polynomials is today a fairly mature subject,
after three decades of activity by many people, building on the pioneering
work of Douady and Hubbard. Given a polynomial $p$ of degree $d\geq 2$,
the most important set is the Julia set $J$ consisting of points $z\in\C$
which have no neighborhood in which the family of iterates forms a normal
family in the sense of Montel. Specifically for polynomials, one can
equivalently start with the set $I$ of points which converge to infinity
under iteration (the {\em escaping} points); the complement $K=\C\sm I$
is known as the {\em filled-in Julia set} and consists of the points with
bounded orbits. Then $J=\partial I=\partial K$. 

The most important case is when $J$ and equivalently $K$ are connected.
Then there is a conformal isomorphism $\phi\colon (\C\sm K)\to(\C\sm
\diskbar)$ which conjugates the dynamics of $P$ on $I=\C\sm K$ to the
dynamics of $z\mapsto z^d$ on $\C\sm\diskbar$. The goal is to show that
the inverse Riemann map $\psi=\phi^{-1}\colon\C\sm\diskbar\to I$ extends
continuously to the boundaries as a continuous surjection
$\psi\colon\partial\disk\to J$; this map would provide a topological
semiconjugacy between the dynamics of $z^d$ on $\partial\disk$ to
the dynamics of $p$ on $J$. The set $I$ is canonically foliated
into {\em dynamic rays} $R_\theta=\psi\left((1,\infty)e^{2\pi
i\theta})\right)$ for $\theta\in\R/\Z$. The dynamic ray $R_\theta$ {\em
lands} at $z\in J$ if the limit $\lim_{r\searrow 1}\psi(re^{2\pi i
\theta})$ exists and equals $z$. The statement that $\psi$ extends
continuously to $\partial\disk$ means that every dynamic ray lands, and
the landing points depend continuously on the angle. By Carath\'eodory's
theorem, this is true if and only if $J$ is locally connected. In this
case, every point $z\in J$, together with its dynamics, is described by
which dynamic rays land at $z$, and this provides a complete description
of the topological properties of the dynamics of $p$ on $J$.

For transcendental entire functions $f$, the set $I$ of escaping points
is equally important as for polynomials, but it is much harder to
understand: the set $I$ is never empty, and it is never a neighborhood of
infinity (because $\infty$ is an essential singularity), so there is no
Riemann map providing convenient coordinates;  but we still have
$J=\partial I$ \cite{Eremenko}. In many cases, $I$ has no interior and
$J=\ovl I$. Eremenko~\cite{Eremenko} has asked whether every (path)
component of $I$ was unbounded. This has been settled in the affirmative
only for the cases of exponential maps $\lambda e^z$ \cite{Escaping} and
for the cosine family $ae^z+be^{-z}$ \cite{Guenter}: in both cases, every
path component of $I$ consists of a single curve which terminates at
$\infty$, and whose other end might or might not land in $\C$; if it does,
the landing point might or might not escape. Only for these rather special
maps is it currently known that the escaping points are organized in the
form of dynamic rays. For a larger class of maps, the existence of some
curves consisting of escaping points was shown in \cite{DT}. A description
of escaping points for a large class of entire function is currently work
in progress by Rottenfu{\ss}er. 

In this paper, we provide the first case of transcendental entire functions
in which it is known that the escaping points are organized in the form of
dynamic rays, such that every dynamic ray lands and every point in $\C$
(which equals the Julia set) is either on a dynamic ray or a landing point
of at least one dynamic ray. To our knowledge, this is the first example
of a transcendental function for which the Julia set equals $\C$ and the
dynamics on all of $\C$ is described in terms of symbolic dynamics.

In addition, we obtain a surprising result about the Hausdorff dimensions
of the involved sets: it turns out that the union of all dynamic rays has
dimension $1$, while each of the remaining points in $\C$ (almost all
points in $\C$ in a very strong sense!) is the landing point of one or
several rays (each ray of course has dimension $1$, and so does their
union!).

In \cite{Bogusia}, Karpi\'nska had shown that for exponential maps
$z\mapsto\lambda\exp(z)$ with attracting fixed points (in particular,
$0<\lambda<1/e\subset\R$), the set of landing points of $R$ has Hausdorff
dimension $2$, while $R$ has dimension $1$. This was extended in
\cite{Escaping} to arbitrary exponential maps: not all rays land, and not
all rays which land have escaping landing points, but the union of all
rays still has dimension $1$, while the set of escaping landing points has
dimension $2$ (but planar measure zero). In \cite{Guenter}, the analogous
result was shown for arbitrary maps $z\mapsto ae^z+be^{-z}$, except that
the escaping landing points of $R$ have even positive measure (using
McMullen's result
\cite{McMullen}). The example given in the present paper is maximal
possible in the plane.
In \cite{Hausdorff}, our results are discussed and illustrated in special
cases, together with an introductory discussion of Hausdorff dimension and
background from complex dynamics.

{\sc Acknowledgement}.
The inspiration for this work occured during a conference on Tenerife in
early 2002, sponsored by the Mittag-Leffler institute in Djursholm/Sweden.
I am most grateful for their support. On this meeting and elsewhere, I
have enjoyed fruitful discussions with Lukas Geyer.

\section{Notation and Background}

Define the sets $\Z_L:=\{\dots,-2_L, -1_L, 0_L, 1_L, 2_L, \dots\}$ and
$\Z_R:=\{\dots,-2_R, -1_R, 0_R, 1_R, 2_R, \dots\}$ (two disjoint copies of
$\Z$) and set $\Sym_0:=(Z_L\cup Z_R)^\N$. To simplify notation, numbers
$s_L$ and $s_R$ (with $s\in\Z$) are treated just like ordinary integers in
arithmetic operations such as $2\pi i s_L=2\pi i s$ or $|s_R|=|s|$, etc..

The space $\Sym_0$ is endowed with the shift map
$\sigma\colon\Sym_0\to\Sym_0$. For every $\s=s_1s_2s_3\dots\in\Sym_0$,
define
\[
t_\s:=\inf\left\{t>0\colon \lim\frac{|s_k|}{F^{\circ k}(t)}=0\right\}
\in\R_0^+\cup\{\infty\}
\,\,.
\]
Define $F\colon\R_0^+\to\R_0^+$ via $F(t):=e^t-1$. A sequence
$\s=s_1s_2s_3\dots\in\Sym_0$ is called {\em exponentially bounded} if
there is an $x\in\R$ such that $|s_k|\leq F^{\circ (k-1)}(x)$ for all $k$.
By \cite[Theorem~4.2 (1)]{Escaping}, a sequence $\s\in\Sym_0$ is
exponentially bounded if and only if $t_\s<\infty$. Let
$\Sym\subset\Sym_0$ be the space of all exponentially bounded sequences.

Define $E(z):=ae^z+be^{-z}$ with $a,b\in\Cstar$ and $I:=\{z\in\C\colon
E^{\circ k}(z)\to\infty \mbox{ as $k\to\infty$}\}$ (the set of escaping
points). 

The following two theorems are the main results in \cite{Guenter}, and
they hold for every map $E(z)=ae^z+be^{-z}$ with $ab\neq 0$ (if one or
both critical orbits escape, then the statements have to be modified
slightly in a natural way).

\begin{theorem}[Existence of dynamic rays]
\label{ThmExistenceRays} \lineclear
Suppose that no critical orbit escapes. Then for every exponentially
bounded $\su\in\Sym$ there exists a unique injective curve 
$g_{\su}:(t_{\su},\infty) \rightarrow I $ consisting of escaping points
such that 
\begin{eqnarray}
g_\su(t)&=&t-\alpha+2\pi is_1+o(1) \qquad \mbox{as $t\to\infty$, if
$s_1\in\Z_R$}
\\
g_\su(t)&=&-t+\beta+2\pi is_1+o(1) \qquad \mbox{as $t\to\infty$, if
$s_1\in\Z_L$}
\\
E(g_{\su}(t))&=&g_{\sigma(\su)}(F(t)) \qquad\qquad\qquad \mbox{ for all $t>
t_{\su}$} \,\,.
\end{eqnarray}
Moreover, for every $t>t_\su$ the orbit of $\ray(t)$ satisfies the
following asymptotics as $k\to\infty$:
\begin{equation}
E^{\circ k}(\ray(t)) =
\left\{ 
\begin{array}{ll}
 F^{\circ k}(t) -\alpha+2\pi i s_{k+1} +o(1)
& \mbox{if $s_{k+1}\in\Z_R$} 
\\
-F^{\circ k}(t) +\beta+2\pi i s_{k+1} +o(1)
&
\mbox{if $s_{k+1}\in\Z_L$ .} 
\end{array}
\right.
\label{EqAsymptotics}
\end{equation}
In particular, the orbit $z_k:=E^{\circ k}(\ray(t))$ satisfies
\begin{equation}
\frac{\log^+|\Im(z_k)|}{\log|\Re(z_k)|} \longrightarrow 0 \,\,.
\label{EqLogEstimate}
\end{equation}
\end{theorem}

\begin{theorem}[Escaping points are organized in rays]
\label{ThmEscapingPointsOnRays} \lineclear
Suppose that no critical orbit escapes. Then for every escaping point $w$ 
there exists a unique exponentially bounded external address $\su$ and a
unique potential $t \geq t_{\su}$ such that exactly one of the following
holds:
\begin{itemize}
\item 
either $t > t_{\su}$ and $w= g_{\su}(t)$,
\item  
or $t=t_{\su}$ and the dynamic ray $g_{\su}$  lands at $w$ such that
$w$ and the ray $g_{\su}$ escape uniformly.
\end{itemize}
In particular, every path component of $I$ in $\C$ is a dynamic ray,
possibly together with the escaping landing point of the ray.
\end{theorem}

We will often need the union of all dynamic rays:
\[
R:=\bigcup_{\s\in\Sym} \ray((t_\s,\infty)) \subset I\,\,.
\]

The following lemma is shown in \cite{Guenter}. 

\begin{lemma}[Horizontal Expansion]
\label{LemHorizExp} \lineclear
For every $a,b\in\Cstar$ and $h>0$, there is an $\eta>0$ with the
following  property: if $(z_k)$ and $(w_k)$ are two orbits under $E$ such
that $|\Im(z_k)-\Im(w_k)|<h$ for all $k$, and if\/
$|\Re(z_1)|-|\Re(w_1)|>\eta$,  then $z_1\in R$:
there is an $\s\in\Sym$ and a $t>t_\s$ such that $z_1=\ray(t)$. 
\end{lemma}

While the proof is technically unpleasant, it has a very simple idea: if
$|\Re(z_1)|>|\Re(w_1)|+\eta$, then the formula $|E(z)|\approx
c\exp(|\Re(z)|)$ (with $c\in\{|a|,|b|\}$) shows that $|z_2|\gg|w_2|$, and
since $z_2$ and $w_2$ have essentially equal imaginary parts, this means
that $|\Re(z_2)|\gg|\Re(w_2)|$. By induction, this shows that the real
parts of $(z_k)$ grow much faster than those of other points with
comparable imaginary parts. This implies that $z_1$ cannot be the landing
point of a dynamic ray: if $z_1=\ray(t_\s)$ for some $\s\in\Sym$, then
$z_1$ escapes much faster than other points on the same ray $\ray$, and
this leads to a contradiction to the asymptotics in
Theorem~\ref{ThmExistenceRays}. Therefore,
Theorem~\ref{ThmEscapingPointsOnRays} implies that $z_1$ is a point on a
dynamic ray.

\section{Landing of Dynamic Rays}
\label{SecLanding}

From now on, we will restrict to the special case of {\em postsingularly 
preperiodic} maps: those maps $z\mapsto E(z):=ae^z+be^{-z}$ for
which the two critical values  $\pm 2\sqrt{ab}$ are strictly
preperiodic. The easiest such maps are $z\mapsto k\pi \sinh(z)$ (with
$k\in\Z\sm\{0\}$) for which the critical values are $\pm k\pi i$, and both
map to the repelling fixed point $0$ (the maps $z\mapsto k\pi\sin(z)$ are
the same maps in a rotated coordinate system). Slightly more generally, if
$a=-b$ is such that $a(1-\sinh(2a))=i\pi k$ with $k\in\Z\sm\{0\}$, then
$E$ has both critical values mapping to fixed points.
Since such maps have no finite asymptotic values, while all critical
values are strictly preperiodic, it is well known that all periodic orbits
are repelling. In particular, $0$ is a repelling fixed point with
$E'(0)=k\pi$. 

In the rest of the paper, we will need  
\[
P:=\bigcup_{k\geq 0}E^{\circ k}\left(\{v, v'\}\right)
\qquad
\mbox{(the finite postsingular set)}
\]
and $V:=\C\sm P$; this carries a unique normalized hyperbolic metric.

\begin{lemma}[Dynamic Rays at Critical Values]
\label{LemDynRayCritValues} \lineclear
For every postsingularly preperiodic map $E$, there are preperiodic 
dynamic rays which land at the two critical values, at least one ray at
each critical value.
\end{lemma}

Note that in the special case of $E(z)=k\pi\sinh(z)$, both $\R^+$ and
$\R^-$ are easily seen to be periodic dynamic rays landing at $0$; since
for these maps, the critical values map to $0$, the claim of the lemma
is obvious.

\proof
Choose some periodic $p_0\in P$ and consider a continuously differentiable
curve $\gamma_0\colon[0,\infty)\to\C\sm P$ with $\gamma_0(0)=p_0$ and
$\gamma_0(t)\to\infty$ as $t\to\infty$. Using this, we construct a family
of curves $\gamma_n\colon[0,\infty)\to\C$ such that
$E(\gamma_{n+1}(t))=\gamma_n(F(t))$ for all $t$; the curve $\gamma_{n+1}$
is uniquely determined by requiring that $p_{n+1}:=\gamma_{n+1}(0)$ is the
unique periodic point with $E(p_{n+1})=p_n$. We claim that there are $k,l$
such that $\gamma_l$ and $\gamma_{l+k}$ are homotopic rel $P$. Observe
that homotopies of $\gamma_n$ lift to homotopies of $\gamma_{n+1}$.
We will use the hyperbolic metric in $V:=\C\sm P$.

\begin{figure}[htbp]
\includegraphics[height=60mm]{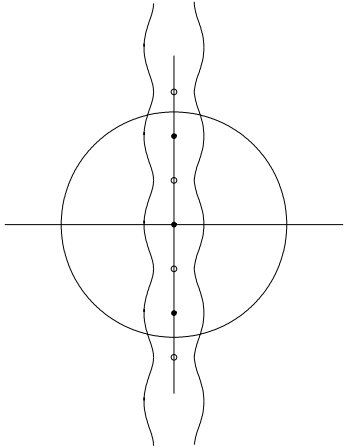}
\caption{Illustration of the proof of Lemma~\ref{LemDynRayCritValues} for
the map $E(z)=\pi\sinh(z)$. Drawn are the circle $C$ and the two preimage
curves $L_1$ and $L_2$. Solid dots indicate the two critical values $\pm
i\pi$ as well as their common image, the repelling fixed point $0$. Some
critical points are indicated by small circles.}
\label{Fig Dyn Rays at Critical Values}
\end{figure}

Let $C$ be a circle of large radius $\rho$, such that $C$ surrounds all of
$P$. Then $E^{-1}(C)$ consists of two unbounded curves $L_1$ and $L_2$, 
one in the right half plane and one in the left half plane, such that both
lines have bounded real parts, while the imaginary parts tend to
$\pm\infty$, and both curves are $2\pi i$-translation invariant. Moreover,
all of $P$ is contained in the unique connected component of
$\C\sm(L_1\cup L_2)$ with bounded real parts. There is
a $\delta>0$ such that every $z\in L_i$ which is surrounded by $C$ can be
connected to $C$ by a curve in $V$ of hyperbolic length of at most
$\delta$; see Figure~\ref{Fig Dyn Rays at Critical Values}.

Let $M$ be the period of $p_0$. Choose $\eps$ small enough so that all
Euclidean disks $D_\eps(p)$ have disjoint closures for all $p\in P$, and
such that $E^{\circ M}(\ovl{D_\eps(p_0)})\supset D_\eps(p_0)$. Let
$\ell_0$ be the hyperbolic length of $\gamma_0$ between $D_\eps(p_0)$ and
$C$. Then the hyperbolic length of $\gamma_1$ between
$E^{-1}(D_\eps(p_0))$ and some $L_i$ is less than $\ell_0$; in fact, there
is an $\eta<1$ such that this hyperbolic length is less than $\eta\ell_0$
because of uniform contraction on compact sets of $V$. After a homotopy of
$\gamma_1$ we may assume that the hyperbolic length of $\gamma_1$ between
$E^{-1}(D_\eps(p_0))$ and $C$ is less than $\eta\ell_0+\delta$, and
$\gamma_1$ intersects $C$ and $\partial E^{-1}(D_\eps(p_0))$ only once
each. After $M$ iterations, it follows that (after an appropriate
homotopy) the hyperbolic length of $\gamma_M$ between $C$ and
$D_\eps(p_0)$ is less than $\eta^M\ell_0+M\delta$. If $\ell_0$ is
sufficiently large, then $\eta^M\ell_0+M\delta<\ell_0$; therefore, the
hyperbolic lengths of all $\gamma_{kM}$ (for $k\in\N$) between
$D_\eps(p_0)$ and $C$ are uniformly bounded above; but this implies that
there are only finitely many homotopy classes available. (For clarity of
exposition, this argument has assumed that every curve intersects $C$ and
$\partial D_\eps(p_0)$ only once; however, this is only a superficial
problem; compare \cite{DynRaysExpo}).

Once we know that we have a curve $\gamma_0$ which is homotopic to
$\gamma_M$ up to homotopy, we stop applying homotopies and consider each
$\gamma_{n+1}$ as a preimage of $\gamma_n$. It then follows that the
curves $\gamma_{kM}$ converge as $k\to\infty$ to a periodic dynamic ray
landing at $p_0$; the details are the same as in
\cite[Section~6]{DynRaysExpo}. The claim follows.
\qed

Using two preperiodic dynamic rays $\ray$ and $\rayp$ landing at the
critical  values $v$ and $v'$, we can introduce a dynamical partition as
follows:  let $U':=\C\sm(\ray\cup\ray'\cup\{v,v'\})$; since $U'$ is simply
connected and contains no singular values, it follows that $U:=E^{-1}(U')$
consists of countably many connected components $W$ such that
$E\colon W\to U'$ is a conformal isomorphism for every $W$. Every
critical point has local mapping degree $2$, so it is the landing point of
exactly two pre-image rays of $\ray$, $\rayp$. More precisely, if $c$ is a
critical point with $E(c)=v$, say, then $c$ is the landing point of two
preimage rays of $\ray$, such that one preimage ray has real parts tending
to $+\infty$ and the other preimage ray has real parts tending to
$-\infty$. The analogous fact is true for the critical points $c'$ with
$E(c')=v'$. The reason is that critical points are spaced at distances
$i\pi$, and every connected component of $U$ must be mapped by $E$ onto
$U'$. Therefore, the dynamical quotient $\C/(2\pi i\Z)$ contains exactly
two components of $U$. The connected components of $U$ can thus be
described with labels $\u\in \Z/2 =\Z\dot\cup(\Z+\frac 1 2)$ such that
$U_{\u+1}$ is the $2\pi i$-translate of $U_\u$ for all $\u\in\Z/2$, and
$\partial U_\u\cap\partial U_{\u+\frac 1 2}\neq\emptyset$ for all
$\u$. This way, we can define {\em itineraries} for every orbit
$(z_k)$ which stays in $U$ entirely: the itinerary of $z_1$ is the sequence
$\u_1\u_2\u_3\dots$  such that $z_k\in U_{\u_k}$ for each $k$. 

\begin{figure}[htbp]
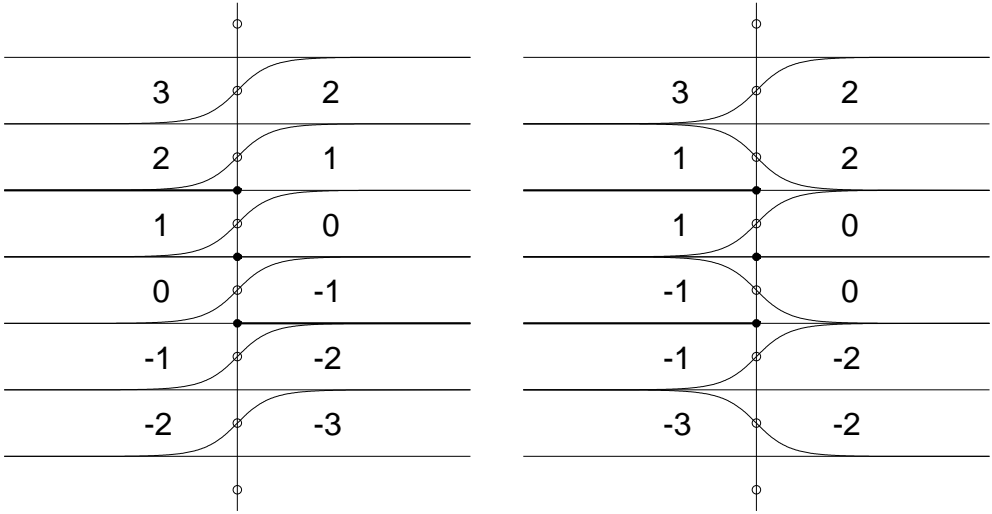

\includegraphics[width=0.49\textwidth]{ray_set1.ps}
\includegraphics[width=0.49\textwidth]{alternating_rays.ps}
\caption{The partition of $\C$ by preimages of rays landing at the
critical values, together with the labels of the components, 
shown for the map $E(z)=\pi\sinh(z)$. Solid dots indicate
again the two critical values $\pm\pi i$ and their common image $0$; small
circles indicate the critical points. Heavy lines indicate one dynamic ray
landing at each of the two critical values (in this simple case, these
rays are horizontal lines); the curved lines are the preimages of these
dynamics rays and land at the critical points. The left and right pictures
show different choices for the used rays at the critical values, and hence
different patterns of the rays landing at the critical points.
}
\label{FigPartition}
\end{figure}

By Theorem~\ref{ThmExistenceRays}, every dynamic ray $\ray$ has
asymptotically constant imaginary parts:
$\lim_{t\to\infty}\Im(\ray(t))$ exists in $\R$. It follows that every
component $U_\u$ has bounded height: there is a number $h>0$ such that for
every $\u\in\Z/2$ and every $z,w\in U_\u$, $|\Im(z)-\Im(w)|<h$.
Therefore, knowing the itinerary of an orbit in $U$ means knowing the
imaginary parts of this orbit up to additive errors of less than $h$.

\begin{lemma}[``Almost'' Every Point in $R$]
\label{LemAlmostInR} \lineclear
For every itinerary $\u_1\u_2\u_3\dots$, among all points $z\in \C$
with orbit in $U$ and with itinerary $\u_1\u_2\u_3\dots$, there are at 
most two which are not in the set $R$ of rays; and if there are two, then
both escape.
\end{lemma}
\proof
Recall that $V:=\C\sm P$ is the complement of the finite postcritical set
and carries a unique normalized hyperbolic metric.
Since  $E^{-1}(V)\subset V$ is a strict inclusion, every branch of $E^{-1}$
(along any curve of finite length) is a strict contraction with respect to
the hyperbolic metric on $V$ (on domain and range). Denote the hyperbolic
distance on $V$ between two points $z,w\in V$ by $d_V(z,w)$. Let $\xi>0$
be such that all $p\in P$ have $|\Re(z)|<\xi$. 

Let $w_1\neq z_1$ be two points with orbits in $U$ and with common
itinerary $\uu:=\u_1\u_2\u_3\dots$ and suppose that both are not in $R$.
Since $E$ restricted to any $U_\u$ is injective, it follows that $z_k\neq
w_k$ for all $k$. Moreover, $|\Im(z_k)-\Im(w_k)|<h$ for every 
$k$, where $h$ bounds the height of every $U_\u$.

By Lemma~\ref{LemHorizExp}, there is a constant $\eta>0$ with the following
property: if there is an index $k$ such that
$\left||\Re(z_k)|-|\Re(w_k)|\rule{0pt}{10pt}\right|>\eta$, then at least
one of $z_k$ and
$w_k$ is in $R$, hence $z_1\in R$ or $w_1\in R$. Therefore, if $z_1$ and
$w_1$ do not belong to $R$, then 
$\left||\Re(z_k)|-|\Re(w_k)|\rule{0pt}{10pt}\right|<\eta$ for every $k$. 
Therefore, either both orbits $z_k$ and $w_k$ are bounded, or both escape,
or both are unbounded without escaping. We treat the three cases
separately.

If both orbits are unbounded but not escaping, then they must infinitely
often visit some compact set $K\subset V$. If $c:=\max\{|a|,|b|\}$, then
$|E(z)|\leq c \exp|\Re(z)|+c$, so there is an $x>\xi+\eta+1$
such that infinitely many $z_k$ satisfy
\[
x\leq|\Re(z_k)|\leq ce^x+c \,\,.
\]
Since always $\left||\Re(z_k)|-|\Re(w_k)|\rule{0pt}{10pt}\right|<\eta$ and
$|\Im(z_k-w_k)|<h$, it follows that for those indices $k$, $d_V(z_k,w_k)$
are uniformly bounded above by some number $d_x>0$. But the pull-back
steps are contracting, uniformly on the compact set $K$, and this implies
that $z_1=w_1$, a contradiction. 

If both orbits $(z_k)$ and $(w_k)$ escape, then the sequences $|\Re(z_k)|$
and $|\Re(w_k)|$ tend to $\infty$ in such a way that
$\left||\Re(z_k)|-|\Re(w_k)|\rule{0pt}{10pt}\right|<\eta$. \linebreak
If
$\left|\Re(z_k)-\Re(w_k)\rule{0pt}{10pt}\right|<\eta$ infinitely often,
then $|z_k-w_k|$ would be bounded above infinitely often while
$|z_k|\to\infty$, and this would imply that $d_V(z_k,w_k)\to 0$ at least
for a subsequence. But as above, we would then have $z_1=w_1$, again a
contradiction. Therefore, for sufficiently large $k$ the real parts of
$z_k$ and $w_k$ always have different signs, and given $z_k$, there is
just one choice for $w_k\neq z_k$. Since the dynamics is injective for the
set of points with identical itineraries, there can be at most two
escaping orbits which are not on rays and which have the same itinerary.

The remaining case is that both orbits are bounded. 
To treat this case, observe first that $P\subset U$: all points in $P$
have well-defined itineraries because both critical values are strictly
preperiodic and the only non-escaping points on the partition boundaries
are the critical points. 

Since both critical orbits land on repelling cycles, it is quite easy to
see that no point in $P$ shares its itinerary with any other non-escaping
point, so none of the two orbits $(z_k)$ and $(w_k)$ can land on $P$.  If
at least one of the two orbits, say $(w_k)$, does not accumulate on $P$,
then there is a compact subset $K\subset V$ with bounded hyperbolic
diameter which contains all $w_k$ and infinitely many $z_k$ (because the
periodic orbits in $P$ are repelling), so this implies again that
$z_1=w_1$. 

The last case is that both orbits accumulate at $P$. Choose some
periodic $p\in P$ in the accumulation set of $(z_k)$, and $\eps>0$ such
that the disk $D_{2\eps}(p)\subset U$ and $D_{2\eps}(p)\cap P=\{p\}$.
There is an $n>0$ such that all points in $D_{\eps}(p)$ have common
itinerary with $p$ for at least $n$ entries. By choosing $\eps>0$
sufficiently small, we may assume that $n$ is large enough so that no two
points in $P$ have common itineraries for $n$ entries.

Let $m$ be the period of $p$
and choose $\eps'>0$ such that $E^{\circ m}(\ovl D_{\eps'}(p))\subset
D_\eps(p)$. Then there are infinitely many $k$ such that $z_k\in
D_\eps(p)\sm D_{\eps'}(p)$. Recall that the orbits $(z_k)$ and $(w_k)$ are
bounded. Since bounded $d_V(z_k,w_k)$ would imply $z_1=w_1$, it follows
that for every subsequence $z_{k_l}\in D_\eps(p)\sm D_{\eps'}(p)$, the
corresponding sequence $w_{k_l}$ must converge to $P$; more precisely,
this sequence must converge to $p$ because this is the only point in $P$
whose itinerary coincides with that of $z_{k_l}$ long enough. Therefore,
if $M_l$ is the number of common entries in the itineraries of $p$ and
$w_{k_l}$ (and hence of $z_{k_l}$), then clearly $M_l\to\infty$. Then
there are numbers $M'_l\leq M_l$ such that $w_{k_l+M'_l}=E^{\circ
M'_l}(w_{k_l})\in D_\eps(p)\sm D_{\eps'}(p)$, and $M'_l\to\infty$.

Similarly as above, the points $z_{k_l+M'_l}$ must converge to $p$; but
this is impossible:  let $U''\subset U$ be the set of points whose
itineraries coincide with that of $p$ for at least $m$ steps; then
$E^{\circ m}\colon U''\to U$ is injective. If we had
$z_{k_l+M'_l}=E^{\circ M'_l}(z_{k_l})\in D_{\eps'}(p)$, then we could pull
back $M'_l$ times, and $z_{k_l}$ would have to be very close to a point in
$P$, which is a contradiction.This excludes the case that $z_k$ and $w_k$
are bounded and proves the claim.
\qed

\begin{theorem}[Dynamic Rays Land]
\label{ThmDynRaysLand} \lineclear
For every postsingularly preperiodic map $E$, every dynamic ray lands in
$\C$.
\end{theorem}
\proof
Let $\ray$ be a dynamic ray, which is a curve
$\ray\colon(t_\s,\infty)\to\C$ (even a $C^\infty$ curve) with bounded
imaginary parts and $\Re(\ray(t))=\pm\infty$ as $t\to\infty$
(Theorem~\ref{ThmExistenceRays}). By \cite[Proposition~6.6]{Guenter},
$t_\s>0$ implies that $\ray$ lands at an escaping point, so we may suppose
that $t_\s=0$.

Let $L_\s\subset\Cbar$ be the limit set of $\ray$: this is the set of
all possible limits of $\ray(t_n)$ as $t_n\searrow 0$. It is well known
that $L_\s=\bigcap_{t>0}\ovl{\ray((0,t))}$, which implies that
$L_\s\subset\Cbar$ is compact and connected.

If $\ray$ is one of the rays bounding the partition $U_\u$, then it lands
at a critical point by definition, and there is nothing to show.
Otherwise, the entire ray is contained in a single domain $U_\u$, hence
$L_\s\subset\ovl{U_\u}$.

Pick any $w\in L_\s$. First we treat the case that $w\in R$, i.e.\ there
are an external address $\s''$ and a $t''>t_{\s''}$ such that
$w=\raypp(t'')$. By (\ref{EqAsymptotics}), we have the asymptotics 
$E^{\circ k}(w)=\pm F^{\circ k}(t'')+2\pi i s''_{k+1}+O(1)$ as
$k\to\infty$. 

Choose any $t'\in(0,t_{\s''})$. Since $\ray(t')$ obeys 
similar asymptotics, there is an $M\in\N$ such that for all $m\geq M$
\[
\left||\Re(E^{\circ m}(w))|-|\Re(E^{\circ
m}(\ray(t')))|\rule{0pt}{11pt}\right|>\eta+1
\]
with the constant $\eta>0$ from Lemma~\ref{LemHorizExp}
(using the height $h$ of the fundamental domains $U_{\u}$). 
But we may choose $t$ sufficiently close to $0$ so that 
$\ray(t)$ is close enough to $w$ such that
\[
|E^{\circ M}(\ray(t))-E^{\circ M}(w)|<1
\,\,,
\]
hence 
\[
\left||\Re(E^{\circ M}(\ray(t)))|
-|\Re(E^{\circ M}(\ray(t')))|\rule{0pt}{11pt}\right|>\eta
\,\,.
\]
By Lemma~\ref{LemHorizExp}, this means that $\ray(t)$ escapes with much
greater real parts than $\ray(t')$; but this contradicts the asymptotics
(\ref{EqAsymptotics}) of the two orbits under the condition $t<t'$.
This excludes the possibility that $w$ is on a dynamic ray, so
$L_\s\subset U_\u\cup C\cup\{\infty\}$, where $C$ denotes the set of
critical points of $E$ (which are the only boundary points of $U_u$ that
are not on rays). 

If the point $w\in L_\s$ does not have a well-defined itinerary, then it
must be either $\infty$ or one of the countably many points on the
backwards orbits of the two critical values. All other points in $L_\s$
have identical itineraries, but none are on dynamic rays, so by
Lemma~\ref{LemAlmostInR}, there can be at most two such points. It follows
that $L_\s$ is countable. However, since $L_\s$ is connected, it contains
at most one point, so $\ray$ lands.

The landing point cannot be $\infty$: otherwise, there would be potentials
$t'>t>0$ such that $|\Re(\ray(t))|-|\Re(\ray(t'))|>\eta$, and again by
Lemma~\ref{LemHorizExp}, this would mean that $\ray(t)$ escapes with much
faster real parts than $\ray(t')$, again a contradiction. Therefore,
$\ray$ lands at some point in $\C$.
\qed

\begin{theorem}[Every Point is Landing Point Or In $R$]
\label{ThmPointLanding} \lineclear
Every $z\in\C$ is either in $R$, or the landing point of at least
one dynamic ray.
\end{theorem}
\proof
We may assume that the itinerary of $z$ is well-defined (or $z$ would
eventually map either onto a dynamic ray on the partition boundary, or
onto the landing point of such a ray). We may also assume that $z$ does
not escape, because every escaping point is either on a ray
or the landing points of a ray (Theorem~\ref{ThmEscapingPointsOnRays}).

For $k\geq 1$, let $V_k\subset\C$ be
the set of points for which at least the first $k$ entries in their
itineraries are well-defined and equal to the itinerary of $z$. Clearly,
$E^{\circ k}\colon V_k\to \C$ is a univalent map with connected unbounded
image, so each $V_k$ is connected and unbounded as well. Moreover, 
$V_{k+1}\subset V_k$ implies that the sets $\ovl V_k\subset\Cbar$ form a
nested sequence of compact and connected sets containing $\{w,\infty\}$,
so $\bigcap_{k\in\N} \ovl V_k$ is also a compact connected set containing
$\{w,\infty\}$. 

Each $V_k$ is bounded by finitely many dynamic rays, say at external
addresses $\s^{k,i}$, together with their landing points on the backwards
orbits of the critical values $\{v,v'\}$. 

If the itinerary of $z$ equals the itinerary of one of the rays
$\rayki$, then the itinerary of $E^{\circ k'}(z)$ equals the itinerary of
a point in $P$ for sufficiently large $k'$, hence $E^{\circ k'}(z)\in P$
and $z$ is the landing point of a dynamic ray. 

We may thus focus on the case that the itinerary of $z$ differs from the
itineraries of all rays $\rayki$; however, the number of common entries in
the itineraries of $\rayki$ and $z$ will tend to $\infty$ as $k\to\infty$.
It now follows that for every ray $\rayki$, there is a $k'>k$ such that
$\rayki\cap\ovl V_{k'}=\emptyset$. For $k\in\N$, let $\s^k$ be an external
address of a dynamic ray in $V_k$.  One can extract a subsequence
from the $\s^k$ which converges pointwise to an external address $\s$.
This address is necessarily exponentially bounded, and the dynamic
ray $\ray$ has a well-defined itinerary which equals that of $z$. 

By Theorem~\ref{ThmDynRaysLand}, the ray $\ray$ lands at some point
$z'\in\C$. Then $z$ and $z'$ have identical itineraries, and both are not
on rays. Since $z$ does not escape, Lemma~\ref{LemAlmostInR} implies that
$z=z'$.
\qed

Set $L:=\C\sm R$: the set of landing points of rays. Since the boundary of
the partition defining itineraries consists of rays landing at critical
points, it follows that every point in $L$ either has a well-defined
itinerary, or is on the backwards orbit of the set $P$.

\section{Dynamics and Dimension}
\label{SecDynamicsDimension}

It is known from \cite[Sec.~7]{Guenter} that if both critical orbits are
strictly preperiodic, then almost every orbit escapes: the set $\C\sm I$
has Lebesgue measure zero. Recall that the set $I$ consists of the rays
together with the landing points of some of the rays, and that
$R$ denotes the union of the rays. Surprisingly, the set $R$ has Hausdorff
dimension $1$ \cite{Guenter}!

\medskip
\noindent
{\sc A dimension paradox.}
As a corollary, we have shown the following surprising result: the set $R$
is a set of Hausdorff dimension $1$ consisting of uncountably many curves
(``rays''), and each ray connects $\infty$ to a well-defined landing
point in $L=\C\sm R$. Conversely, every $z\in L$ has one or even several
rays in $R$ connecting $z$ to $\infty$. One might have expected that the
set of landing points should be ``smaller'' than the set of entire rays,
but the opposite is the case: not only does $\C\sm R$ have greater
Hausdorff dimension than $R$, and not only has it positive or even full
Lebesgue measure in $\C$: it is the complement of the 1-dimensional set
$R$\,! In other words, we have a partition of $\C$ into an uncountable
union of disjoint sets $Y_i=\{z_i\}\cup R_i^1\cup R_i^2\cup\dots$, where
each component $Y_i$ consists of one point $z_i\in L$ and one or several
rays $R_i^j$ landing at $z_i$. Every ray $R_i^j$ is a one-dimensional
curve, and the uncountable union $R=\bigcup_{i,j}R_i^j$ still has
dimension $1$, but the union $L=\bigcup_i\{z_i\}=\C\sm R$ is
two-dimensional and has full planar Lebesgue measure as the complement of
$R$.

\medskip
\noindent
{\sc Remark on possible extensions.}
It seems quite likely that for maps $z\mapsto ae^z+be^{-z}$ in which one
or both critical orbits are allowed to be periodic, rather than
preperiodic, similar results hold as for our maps: 
every dynamic ray lands, and every point in the Julia set is the
landing point of one or several dynamic rays. However, because of the
superattracting cycles the Fatou set would be non-empty and the set of
landing points would no longer have full measure (but still positive).
Every point would be either in $R$, a landing point of finitely
many rays, or in a superattracting basin. One might wonder whether the
precise description of the dynamics would allow for analogs of ``pinched
disk'' models \cite{DoCompacts} or some kinds of combinatorial or even
topological renormalization (although non-continuity of any possible
renormalization for exponential maps has recently been shown by
Rempe~\cite{LasseTopology}, and this argument is likely to apply in our
cases too). 

While we have shown that all dynamic rays land, we have not discussed
whether the landing points depend continuously on the external address. In
fact, this depends on the topology used in $\Sym$: as a totally ordered
space, $\Sym$ possesses its ``order topology'' generated by intervals;
with this topology, the landing points do {\em not} depend continuously on
the address because very nearby addresses can still have very different
minimal potentials. It is therefore required to endow $\Sym$ with a
topology which depends on more than finitely many first entries; this is a
matter of separate discussion.

\end{document}